\documentclass{amsart}
\usepackage{amssymb, setspace}

\newtheorem{theorem}{Theorem}[section]
\newtheorem{lemma}[theorem]{Lemma}

\theoremstyle{definition}
\newtheorem{definition}[theorem]{Definition}
\newtheorem{example}[theorem]{Example}

\theoremstyle{remark}
\newtheorem{remark}[theorem]{Remark}

\numberwithin{equation}{section}
\newcommand{\aut}{\mathrm{Aut}}

\newcommand{\Gal}{\mathrm{Gal}}
\newcommand{\ZZ}{\mathbb{Z}}
\newcommand{\NN}{\mathbb{N}}
\newcommand{\QQ}{\mathbb{Q}}
\newcommand{\RR}{\mathbb{R}}
\newcommand{\FF}{\mathbb{F}}
\newcommand{\CO}{\mathcal{O}}

\newcommand{\Hgt}{\mathrm{Height}}

\newcommand{\gen}[2]{\langle{#1},{#2}\rangle}
\begin{document}
\markboth{H.-C. Li}{Height of Some Automorphisms of Local Fields}
\title{Height of Some Automorphisms of Local Fields}

\author{Hua-Chieh Li}

\address{Department of Mathematics, National Taiwan Normal University, Taipei, Taiwan, R.O.C.}
\email{li@math.ntnu.edu.tw}

\begin{abstract}
In this note, we determine which automorphism subgroups of $\mathrm{Aut}_{\mathbb{F}_q}(\mathbb{F}_q((x)))$ are corresponding to $\mathbb{Z}_p$-extensions or $\mathbb{Z}_p\times\mathbb{Z}_p$-extensions of characteristic $0$ fields.
\end{abstract}

\maketitle

%\begin{history}
%\received{(Day Month Year)}
%\accepted{(Day Month Year)}
%\end{history}

\date{}

%    "Communicated by" -- provide editor's name; required.
%\commby{}

%    Abstract is required.

%\keywords{field of norms functor; upper ramification subgroups; lower ramification subgroups, Hasse-Herbrand function.}
%\ccode{Mathematics Subject Classification 2010: 11S15, 11S31}
%\baselineskip=20pt
\section{Introduction}
Let $k=\FF_q$ be a finite field of characteristic $p>0$ and let $L/K$ be a totally ramified abelian extension, where $K$ is a local field with residue field $k$.
Then $G=\Gal(L/K)$ has a decreasing filtration by the upper ramification subgroups $G(r)$, defined for nonnegative $r\in\RR$ (see \cite[IV]{serre79}). Since $G$ is abelian, $L/K$ is arithmetically profinite (see \cite{win83}). This means that for every $r\ge 0$ the upper ramification group $G(r)$ has finite index in $G$. This allows us to define the Hasse-Herbrand function $\psi:\RR^+\to\RR^+$ where $\psi_{L/F}(r)=\int_0^r[G:G(t)]dt$ and $\phi_{L/K}(r)=\psi^{-1}_{L/K}(r)$. The ramification subgroups of $G$ with the lower numbering are defined by $G[r]=G(\phi_{L/K}(r))$.

Let $\aut_k(k((x)))$ denote the group of continuous automorphisms of $k((x))$ which induce the identity map on $k$. A closed abelian subgroup $G$ of $\aut_k(k((x)))$ also has a ramification filtration. The lower ramification subgroups of $G$ are defined by \[G[r]=\{\sigma\in G: v_x(\sigma(x)-x)\ge r+1\}\] for $r\ge 0$. Since $G[r]$ has finite index in $G$ for every $r\ge 0$, the function $\phi_G:\RR^+\to\RR^+$ where $\phi_G(r)=\int_0^r[G:G[t]]^{-1}dt$ is strictly increasing. We define the ramification subgroups of $G$ with the upper numbering by $G(r)=G[\phi_G^{-1}(r)]$.

Wintenberger \cite{win80} has shown that the field of norms functor induces an equivalence between a category whose objects are totally ramified abelian $p$-adic Lie extensions $L/K$, where $K$ is a local field with residue field $k$, and a category whose objects are pairs $(\mathcal{K},G)$, where $\mathcal{K}\simeq k((x))$ and $G$ is an abelian $p$-adic Lie subgroup of $\aut_k(\mathcal{K})$. In short, if $G$ is an abelian $p$-adic Lie subgroup of $\aut_k(k((x)))$, then there is an abelian $p$-adic Lie extensions $L/K$ corresponding to $(k((x)),G)$ by the equivalence of categories given by the field of norms functor. Moreover, the canonical isomorphism from $\Gal(L/K)$ onto $G$ preserves the ramification filtration \cite{laubie07, win83}. This equivalence has been extended to allow $\Gal(L/K)$ and $G$ to be arbitrary abelian pro-$p$ groups by Keating \cite{keat09}. In the following, we will simply say that $G$ is corresponding to $L/K$ if the extension $L/K$ corresponds to $(k((x)),G)$ by the equivalence of categories given by the field of norms functor.

For $\sigma\in \aut_k( k((x)))$, we let \[i(\sigma)=v_x\left(\frac{\sigma(x)}{x}-1\right).\] Moreover, if $\sigma(x)\equiv x\pmod{x^2}$, then we denote $i_n(\sigma)=i(\sigma^{p^n})$.
When $\sigma\in\aut_k(k((x)))$ has infinite order, the sequence $\{i_n(\sigma)\}$ is strictly increasing and attracts many attentions.
In \cite{sen69} Sen proved that for every $n\in\NN$, $i_{n+1}(\sigma)\equiv i_n(\sigma)\pmod{p^{n+1}}$. In \cite{Keating} Keating determines upper bounds for the $i_n(\sigma)$ in some cases and in \cite{LauSai, LauSai2} the authors improve Keating's results using Wintenberger's theory of field of norms \cite{win80, win83}. These results are base on the fact that the automorphism subgroups correspond to $\ZZ_p$-extensions of characteristic $0$ fields in \cite{Keating, LauSai} and correspond to $\ZZ_p\times\ZZ_p$-extensions of characteristic $0$ fields in \cite{LauSai2}. In this note, we determine which automorphism subgroups of $\aut_k(k((x)))$ are corresponding to $\ZZ_p$-extensions or $\ZZ_p\times\ZZ_p$-extensions of characteristic $0$ fields. In the following, we will simply say that an extension $L/K$ is of characteristic $0$ if the characteristic of $K$ is $0$. Likewise, if the characteristic of $K$ is $p$, then we say that the extension $L/K$ is of characteristic $p$.

Motivated by the definition of height of a formal group and height of a $p$-adic dynamical system \cite{li}, we have the following definition.
\begin{definition}
Let $\sigma\in\aut_k(k((t)))$ with $\sigma\equiv x\pmod{x^2}$. We say that the {\em height} of $\sigma$ exists if $\lim_{n\to\infty}i_{n}(\sigma)/i_{n-1}(\sigma)$ is finite and denote by \[\mathrm{Height}(\sigma)=\lim_{n\to\infty}\log_p \frac{i_n(\sigma)}{i_{n-1}(\sigma)}.\]
\end{definition}
Let $G$ be a closed subgroup of $\aut_k(k((x)))$. Our main result shows that if $G$ is isomorphic to $\ZZ_p$ then $G$ corresponds to a characteristic $0$ field extension if and only if every nonidentity element of $G$ has height $1$ and if $G$ is isomorphic to $\ZZ_p\times\ZZ_p$ then $G$ corresponds to a characteristic $0$ field extension if and only if every nonidentity element of $G$ has height $2$.

The proof of our result is based on the following straightforward consequence of Theorem 4 of \cite{mar}.
\begin{lemma}\label{lemmar}
Let $L/K$ be an abelian extension and let $G$ denote the Galois group $\Gal(L/K)$.
\begin{enumerate}
\item If $K$ is of characteristic $p$, then the mapping $\sigma\to\sigma^p$ maps $G(n)$ into $G(pn)$, for all $n\in\NN$.
\item If $K$ is of characteristic $0$ with absolute ramification index $e$, then the mapping $\sigma\to\sigma^p$ induces a homomorphism which maps $G(n)/G(n+1)$ onto $G(n+e)/G(n+e+1)$, for all $n$ large enough.
\end{enumerate}
\end{lemma}

We remark that since $G$ is abelian, every upper ramification break $u$ (i.e. $G(u)\supsetneq G(u+\epsilon)$, $\forall\,\epsilon>0$) is an integer (see for instance  \cite[V]{serre79}). Therefore, we can apply Lemma \ref{lemmar} to the case where $n$ is an upper ramification break of $G$. Moreover, if $K$ is of characteristic $0$ and $\Gal(L/K)$ is a pro-$p$ group, then Lemma \ref{lemmar} (2) shows that the mapping $\sigma\to\sigma^p$ maps $G(n)$ onto $G(n+e)$, for $n$ sufficiently large. Therefore, in this case, if there is no nontrivial $p$-torsion element in $G$, then the mapping $\sigma\to\sigma^p$ induces an isomorphism between $G(n)/G(n+1)$ and $G(n+e)/G(n+e+1)$, for all $n$ large enough. In particular, if $n$ is large enough and $n$ is an upper ramification break of $G$, then $n+e$ is also an upper ramification break of $G$.

\section{$\ZZ_p$-extensions}
Given $\sigma\in\aut_k(k((t)))$ with $\sigma\equiv x\pmod{x^2}$, write $\lim_{n\to \infty} \left(i_n(\sigma)/p^n\right)
= \left(p/(p-1)\right) e.$ It is well-known that either $e$ is a
positive integer or $e =\infty$ (see for instance \cite{win04}). Moreover, $e$ is a
positive integer if and only if the field extension $E/F$
corresponding to the closed subgroup generated by $\sigma$ is of characteristic 0. In fact, in this case, $e$ is
the absolute ramification index of $F$. If $e$ is finite, then it's clear that
$\lim_{n\to\infty} \left(i_n(\sigma)/i_{n-1}(\sigma)\right) = p$. That
is $\mathrm{Height}(\sigma)=1$. In this section, we will show that the converse is also true.
Thus, for the case $\sigma\in \aut_k(k((x)))$ with $\Hgt(\sigma)=1$, the closed cyclic group generated by $\sigma$ corresponds to a $\ZZ_p$-extension of characteristic $0$ field.

We prove this by contradiction. Suppose that the corresponding $\ZZ_p$-extension is of characteristic $p$. Then it is also true that the $\ZZ_p$-extension corresponding to the closed subgroup $H$ generated by $\sigma^{p^n}$ is of characteristic $p$.
By considering the ramification groups of $H$, we have $\sigma^{p^n}\in H[i_n(\sigma)]\setminus H[i_n(\sigma)+\epsilon]$ and $\sigma^{p^{n+1}}\in H[i_{n+1}(\sigma)]\setminus H[i_{n+1}(\sigma)+\epsilon]$, $\forall\epsilon>0$. Therefore $\phi_H(i_n(\sigma))$ and $\phi_H(i_{n+1}(\sigma))$ are upper ramification breaks of $H$ and hence we can apply Lemma \ref{lemmar} (1) to get $\sigma^{p^{n+1}}\in H(p\,\phi_H(i_n(\sigma)))$. In other words,
\[\phi_{H}(i_{n+1}(\sigma))=i_{n}(\sigma)+\frac{i_{n+1}(\sigma)-i_{n}(\sigma)}{p}\ge p\,\phi_H(i_n(\sigma))=p\,i_{n}(\sigma),\forall\,n\in\NN.\] This says that \[i_{n+1}(\sigma)\ge (p^2-p+1)i_{n}(\sigma),\forall\,n\in\NN,\] and hence contradicts to the assumption that $\lim_{n\to\infty}\frac{i_n(\sigma)}{i_{n-1}(\sigma)}=p.$

Conversely, suppose that $G$ is corresponding to a characteristic $0$ field extension $E/F$ with $e = v_F(p)$ being
the absolute ramification index of $F$. By the definition of lower ramification group, for every $n\in\NN$, $\sigma^{p^n}\in G[i_n(\sigma)]\setminus G[i_n(\sigma)+\epsilon]$ and $\sigma^{p^{n+1}}\in G[i_{n+1}(\sigma)]\setminus G[i_{n+1}(\sigma)+\epsilon]$, for every $\epsilon>0$. On the other hand, by Lemma \ref{lemmar} (2) and the remark following it, when $n$ is large enough if we let $u=\phi_G(i_n(\sigma))$, then $u+e$ is an upper ramification break of $G$. Moreover, since $\sigma^{p^n}\in G[i_n(\sigma)]=G(u)$, $\sigma^{p^{n+1}}\in G(u)^p=\{g^p:g\in G(u)\}=G(u+e).$ In other words, $u+e=\phi_G(i_n(\sigma))+e=\phi_G(i_{n+1}(\sigma)),$ and hence \[e=\phi_G(i_{n+1}(\sigma))-\phi_G(i_n(\sigma))=\frac{1}{p^{n+1}}(i_{n+1}(\sigma)-i_n(\sigma))\] because $G$ is isomorphic to $\ZZ_p$. This is true for all $n$ large enough. Therefore, we conclude that there exists $m\in\NN$ such that for all $n>m$,
\begin{eqnarray*}
  i_n(\sigma) &=& i_m(\sigma)+\sum_{j=m+1}^n(i_{j}(\sigma)-i_{j-1}(\sigma)) \\
   &=& i_m(\sigma)+e(p^{m+1}+\cdots+p^n) \\
   &=& i_m(\sigma)+\frac{ep}{p-1}(p^{n}-p^{m}).
\end{eqnarray*} This shows
\begin{equation}\label{ram1}
  \lim_{n\to\infty}\frac{i_n(\sigma)}{p^n}=\frac{ep}{p-1}.
\end{equation}

We summarize this result as the following.

\begin{theorem}\label{height1}
Suppose that $G\subseteq\aut_k(k((x)))$ is a closed subgroup generated by $\sigma$ which is isomorphic to $\ZZ_p$. Then the following are equivalent:
\begin{enumerate}
\item $\lim_{n\to\infty}\frac{i_{n}(\sigma)}{i_{n-1}(\sigma)}=p$
\item The sequence $\{\frac{i_n(\sigma)}{p^n}\}_n$ converges.
\item $\frac{i_{n+1}(\sigma)-i_n(\sigma)}{i_n(\sigma)-i_{n-1}(\sigma)}=p$ for all $n$ sufficiently large.
\item The $\ZZ_p$-extension corresponding to $G$ is of characteristic $0$.
\end{enumerate}

\end{theorem}

\begin{remark}
Theorem \ref{height1} remains true if we replace $\sigma$ in the statements (1), (2) and (3) by any nonidentity element $\tau\in G$. This is because the closed subgroup of $G$ generated by $\tau$ is a finite index subgroup. In other words, we shows that the $\ZZ_p$-extension corresponding to $G$ is of characteristic $0$ if and only if every nonidentity element $\tau\in G$ has $\Hgt(\tau)=1$.
\end{remark}

\section{$\ZZ_p\times\ZZ_p$-extensions}
In this section we extend the result of the previous section to the case that $G\subseteq\aut_k(k((x)))$ is isomorphic to $\ZZ_p\times \ZZ_p$.
In this case we show that every nonidentity element of $G$ has height $2$ if and only if the $\ZZ_p\times \ZZ_p$-extension corresponding to $G$ is of characteristic $0$.

Let $G$ be a closed subgroup of $\aut_k(k((x)))$ which is isomorphic to $\ZZ_p\times \ZZ_p$ and suppose that for every nonidentity element $\sigma\in G$ we have $\lim_{n\to\infty}i_{n}(\sigma)/i_{n-1}(\sigma)=p^2$. Again, we use method of contradiction to show that the $\ZZ_p\times \ZZ_p$-extension corresponding to $G$ is of characteristic $0$. First, suppose that the corresponding $\ZZ_p\times \ZZ_p$-extension is of characteristic $p$. Then for any two linearly independent elements $\sigma,\tau\in G$, since $\gen{\sigma}{\tau}$ is a finite index subgroup of $G$ (we use $\gen{\sigma}{\tau}$ to denote the closed subgroup of $G$ generated by $\sigma$ and $\tau$), the field extension corresponding to $\gen{\sigma}{\tau}$ is also a characteristic $p$ field extension.
Similarly, for $m,n\in\NN$, the $\ZZ_p\times\ZZ_p$ extension corresponding to $\gen{\sigma^{p^n}}{\tau^{p^m}}$ is also of characteristic $p$.
We consider several cases.

For a given nonidentity $\sigma\in G$, we first suppose that for every $N\in\NN$, there exist $n,m>N$
and $\tau$ of $G$ such that $i_n(\sigma)< i_m(\tau)<i_{m+1}(\tau)\le i_{n+1}(\sigma)$. Notice that the field extension corresponding to the closed subgroup $H=\gen{\sigma^{p^n}}{\tau^{p^m}}$ is also of characteristic $p$. By considering the lower ramification subgroups of $H$, we have
\begin{align*}
  H[i_n(\sigma)]=\gen{\sigma^{p^n}}{\tau^{p^m}}&\supsetneq H[i_n(\sigma)+1]=\cdots =H[i_m(\tau)]=\gen{\sigma^{p^{n+1}}}{\tau^{p^m}} \\
  &\supsetneq H[i_m(\tau)+1]=\cdots =H[i_{m+1}(\tau)]=\gen{\sigma^{p^{n+1}}}{\tau^{p^{m+1}}}.
\end{align*}
Therefore, \[\phi_H(i_m(\tau))=i_n(\sigma)+\frac{i_m(\tau)-i_n(\sigma)}{p}\] and
\[\phi_H(i_{m+1}(\tau))=i_n(\sigma)+\frac{i_m(\tau)-i_n(\sigma)}{p}+\frac{i_{m+1}(\tau)-i_m(\tau)}{p^2}.\]
Since $\tau^{p^m}\in H[i_m(\tau)]=H(\phi_H(i_m(\tau)))$ and
\[\tau^{p^{m+1}}\in H(\phi_H(i_{m+1}(\tau)))\setminus H(\phi_H(i_{m+1}(\tau))+\epsilon),\,\, \forall\epsilon>0,\] Lemma \ref{lemmar} (1) says    \[i_n(\sigma)+\frac{i_m(\tau)-i_n(\sigma)}{p}+\frac{i_{m+1}(\tau)-i_m(\tau)}{p^2}\ge p(i_n(\sigma)+\frac{i_m(\tau)-i_n(\sigma)}{p}).\]
Therefore by $i_{n+1}(\sigma)\ge i_{m+1}(\tau)$ and $i_m(\tau)>i_n(\sigma)$, we have \[i_{n+1}(\sigma)\ge(p^2-p+1)i_m(\tau)+(p^3-2p^2+p)i_n(\sigma)>(p^3-p^2+1)i_n(\sigma).\]
Since for every $N\in\NN$, this is true for some $n>N$, it contradicts to the assumption that $\lim_{n\to\infty}\frac{i_n(\sigma)}{i_{n-1}(\sigma)}=p^2.$

Now suppose that for every $N\in\NN$, there exist $n,m>N$
and $\tau$ of $G$ such that $i_n(\sigma)< i_m(\tau)< i_{n+1}(\sigma)<i_{n+2}(\sigma)\le i_{m+1}(\tau)$. Notice that the field extension corresponding to the closed subgroup $H=\gen{\sigma^{p^n}}{\tau^{p^m}}$ is also of characteristic $p$. By considering the lower ramification subgroups of $H$, we have
\begin{align*}
  H[i_n(\sigma)]=\gen{\sigma^{p^n}}{\tau^{p^m}}&\supsetneq H[i_n(\sigma)+1]=\cdots =H[i_m(\tau)]=\gen{\sigma^{p^{n+1}}}{\tau^{p^m}} \\
  &\supsetneq H[i_m(\tau)+1]=\cdots =H[i_{n+1}(\sigma)]=\gen{\sigma^{p^{n+1}}}{\tau^{p^{m+1}}}\\
  &\supsetneq H[i_{n+1}(\sigma)+1]=\cdots =H[i_{n+2}(\sigma)]=\gen{\sigma^{p^{n+2}}}{\tau^{p^{m+1}}}.
\end{align*}
Therefore \[\phi_H(i_{n+2}(\sigma))=\phi_H(i_{n+1}(\sigma))+\frac{i_{n+2}(\sigma)-i_{n+1}(\sigma)}{p^3}.\]
By Lemma \ref{lemmar} (1), $\phi_H(i_{n+1}(\sigma))\ge p\, \phi_H(i_{n}(\sigma))=p\,i_n(\sigma)$ and since $\sigma^{p^{n+2}}\not\in H[i_{n+2}(\sigma)+\epsilon],$ $\forall\,\epsilon>0$, we have $\phi_H(i_{n+2}(\sigma))\ge p\, \phi_H(i_{n+1}(\sigma))$.
This implies \[\frac{i_{n+2}(\sigma)-i_{n+1}(\sigma)}{p^3}\ge (p-1)\phi_H(i_{n+1}(\sigma))\ge (p-1)p\,i_n(\sigma),\] and hence \[i_{n+2}(\sigma)\ge(p^5-p^4)i_n(\sigma)+i_{n+1}(\sigma).\]
Since for every $N\in\NN$, this is true for some $n>N$, it contradicts to the assumption that $\lim_{n\to\infty}\frac{i_{n+1}(\sigma)}{i_{n}(\sigma)}=p^2.$

Now we only have the following two cases to consider: \begin{enumerate}
\item There exists $N$ such that there is neither $m,n>N$ nor any nonidentity $\tau\in G$ such that $i_{n}(\sigma)<i_{m}(\tau)<i_{n+1}(\sigma)$.
\item There exists $m,n\in\NN$ and a nonidentity $\tau\in G$ such that $i_{n+j}(\sigma)<i_{m+j}(\tau)<i_{n+j+1}(\sigma)$, for all $j\in\NN$.
\end{enumerate}
For the case (1), there exists $N\in\NN$ such that \[G[i_n(\sigma)]\supsetneq G[i_n(\sigma)+1]=\cdots=G[i_{n+1}(\sigma)],\,\forall\,n>N.\] Now let $H=G[i_n(\sigma)]$. Then by the contrapositive assumption,
the field extension corresponding to $H$ is also of characteristic $p$.
Since $\phi_H(i_{n+1}(\sigma))=i_n(\sigma)+(1/p^2)(i_{n+1}(\sigma)-i_n(\sigma))$, again by Lemma \ref{lemmar} (1) we have $i_{n}(\sigma)+(1/p^2)(i_{n+1}(\sigma)-i_{n}(\sigma))\ge pi_{n}(\sigma)$ and hence \[i_{n+1}(\sigma)\ge (p^3-p^2+1)i_{n}(\sigma)\ge(p^2+1)i_{n}(\sigma).\]
This is true for all $n>N$, and hence it contradicts to the assumption that $\Hgt(\sigma)=2.$

For the case (2), for every $j\in\NN$, let $H=\gen{\sigma^{p^{n+j}}}{\tau^{p^{m+j}}}$ and by considering the ramification subgroups of $H$, we have
\[\phi_H(i_{n+j+1}(\sigma))=i_{n+j}(\sigma)+\frac{i_{m+j}(\tau)-i_{n+j}(\sigma)}{p}+\frac{i_{n+j+1}(\sigma)-i_{m+j}(\tau)}{p^2}.\]
Again, by the contrapositive assumption,
the field extension corresponding to $H$ is of characteristic $p$, and hence by Lemma \ref{lemmar} (1)
\[\phi_H(i_{n+j+1}(\sigma))\ge p\phi_H(i_{n+j}(\sigma))=p\,i_{n+j}(\sigma).\] This says that
\begin{equation}\label{eq1}
i_{n+j+1}(\sigma)\ge (p^3-p^2+p)i_{n+j}(\sigma)-(p-1)i_{m+j}(\tau).
\end{equation}
Since $\lim_{j\to\infty}i_{n+j+1}(\sigma)/i_{n+j}(\sigma)=p^2$, for every $1>\epsilon >0$, there exists $j$ large enough such that $p^2-\epsilon<i_{n+j+1}(\sigma)/i_{n+j}(\sigma)<p^2+\epsilon$. Similarly, $p^2-\epsilon<i_{m+j+1}(\tau)/i_{m+j}(\tau)<p^2+\epsilon$. Hence, we can have either $i_{m+j}(\tau)< (p+\epsilon)i_{n+j}(\sigma)$ or $i_{n+j+1}(\sigma)< (p+\epsilon)i_{m+j}(\tau)$. Otherwise $i_{m+j}(\tau)\ge (p+\epsilon)i_{n+j}(\sigma)$ and $i_{n+j+1}(\sigma)\ge (p+\epsilon)i_{m+j}(\tau)$ imply $i_{n+j+1}(\sigma)\ge (p+\epsilon)^2 i_{n+j}(\sigma)>(p^2+\epsilon)i_{n+j}(\sigma)$. Without lose of generality (switching $\sigma$ and $\tau$ if necessary), for every $N\in\NN$ and $1>\epsilon>0$, we can find $j>N$ such that $i_{m+j}(\tau)<(p+\epsilon)i_{n+j}(\sigma)$ and hence by Equation (\ref{eq1}), we get
\[i_{n+j+1}(\sigma)\ge (p^3-p^2+p)i_{n+j}(\sigma)-(p-1)(p+\epsilon)i_{n+j}(\sigma).\] Thus \[i_{n+j+1}(\sigma)\ge (p^3-2p^2+(2-\epsilon)p+\epsilon)i_{n+j}(\sigma).\]
This contradicts to the assumption that $\lim_{j\to\infty}\frac{i_{n+j+1}(\sigma)}{i_{n+j}(\sigma)}=p^2$, for $p\ge 3$.

For the case $p=2$, considering the ramification subgroups
\[
  H[i_{n+j}(\sigma)]\supsetneq H[i_{m+j}(\tau)]\supsetneq H[i_{n+j+1}(\sigma)]\supsetneq H[i_{m+j+1}(\tau)]\supsetneq H[i_{n+j+2}(\sigma)],
\] we have \begin{multline*}
\phi_H(i_{n+j+2}(\sigma))= i_{n+j}(\sigma)+\frac{i_{m+j}(\tau)-i_{n+j}(\sigma)}{2}+\frac{i_{n+j+1}(\sigma)-i_{m+j}(\tau)}{4})\\
             +\frac{i_{m+j+1}(\tau)-i_{n+j+1}(\sigma)}{8}+\frac{i_{n+j+2}(\sigma)-i_{m+j+1}(\tau)}{16}.
           \end{multline*}
Again by the assumption that the corresponding field extension is of characteristic $2$, we have $\phi_H(i_{n+j+2}(\sigma))\ge 2\phi_H(i_{n+j+1}(\sigma))$
and deduce that
\[i_{n+j+2}(\sigma)\ge 8i_{n+j}(\sigma)+4i_{m+j}(\tau)+6i_{n+j+1}(\sigma)-i_{m+j+1}(\tau).\] Again, without lose of generality, for every $N\in\NN$ and $1>\epsilon>0$, we can assume there exists $j>N$ such that $i_{n+j+1}(\sigma)> (4-\epsilon)i_{n+j}(\sigma)$, $i_{m+j+1}(\tau)<(4+\epsilon)i_{m+j}(\tau)$ and $i_{m+j}(\tau)<(2+\epsilon)i_{n+j}(\sigma)$. Therefore, by using $i_{m+j}(\tau)>i_{n+j}(\sigma)$, we get
\[i_{n+j+2}(\sigma)>(28-12\epsilon-\epsilon^2)i_{n+j}(\sigma).\] This contradicts to the assumption that $\lim_{j\to\infty}\frac{i_{n+j+2}(\sigma)}{i_{n+j}(\sigma)}=2^4$. We complete the proof of showing that if every nonidentity element of $G$ is of height $2$, then the $\ZZ_p\times\ZZ_p$-extension corresponding to $G$ is of characteristic $0$.

Conversely, suppose the field extension $E/F$
corresponding to $G$ is of characteristic 0 with $e = v_F(p)$ being
the absolute ramification index of $F$. Then since there is no $p$-torsion element in $G$, by Lemma \ref{lemmar} (2) and the remark following it, there exists an $N$ such that the raise to $p$-th power map $G(u)/G(u+\epsilon)\mapsto G(u+e)/G(u+e+\epsilon)$ is an isomorphism for all $u>N$. In other words, there exists $N\in\NN$ such that for every upper ramification break $u>N$, $[G(u):G(u+\epsilon)]$ is either always $2$ or always $1$. For simplicity, we call the former {\em depth} $2$ case and the latter {\em depth} $1$ case.

For depth $2$ case, it means that
for every $\sigma,\tau\in G$, there exists $n,m\in\NN$ such that $i_{n+j}(\sigma)=i_{m+j}(\tau)$ for all $j\in\NN$.
Therefore, for every $\sigma\in G$, we choose another $\tau\in G$ so that $G[i_n(\sigma)]=G(u_1)=\gen{\sigma^{p^n}}{\tau^{p^m}}$, $G[i_{n+1}(\sigma)]=G(u_2)=\gen{\sigma^{p^{n+1}}}{\tau^{p^{m+1}}}$, where $u_1,u_2$ are upper ramification breaks and $G(u)^p=G(u+e)$ for all $u\ge u_1$.
Let $G'$ be the closed subgroup of $G$ generated by $\sigma$ and $\tau$. It is clear that $G'$ is of finite index over $G$ and hence the field extension $E/F'$ corresponding to $G'$ is also of characteristic $0$. Let $e'$ be the absolute ramification index of $F'$. Since $G'[i]=G[i]\cap G'$, we get $\phi_{G'}(i_{n+1}(\sigma))-\phi_{G'}(i_{n}(\sigma))=(i_{n+1}(\sigma)-i_n(\sigma))/p^{n+m+2}=e'$.
Inductively, we have
\begin{equation}\label{dep2}
  \frac{i_{n+j}(\sigma)-i_{n+j-1}(\sigma)}{p^{n+m+2j}}=e'.
\end{equation}
This shows
\[\frac{i_{n+1}(\sigma)-i_{n}(\sigma)}{i_{n}(\sigma)-i_{n-1}(\sigma)}=p^2\] for all $n$ large enough.
Moreover, since \[i_{n+j}(\sigma)=i_n(\sigma)+(i_{n+1}(\sigma)-i_n(\sigma))+\cdots+(i_{n+j}(\sigma)-i_{n+j-1}(\sigma)),\] we have
\[i_{n+j}(\sigma)=i_n(\sigma)+\frac{p^{2+m+n}e'}{p^2-1}(p^{2j}-1)\]
and hence the limit $\lim_{n\to\infty}\frac{i_n(\sigma)}{p^{2n}}$ exists.

For depth $1$ case, it means that
for every $\sigma\in G$, there exists $\tau\in G$ and $n,m\in\NN$ such that $i_{n+j}(\sigma)<i_{m+j}(\tau)<i_{n+j+1}(\sigma)$ for all $j\in\NN$.
Therefore, for every $\sigma\in G$, we choose another $\tau\in G$ satisfying this condition so that $G[i_n(\sigma)]=G(u_1)=\gen{\sigma^{p^n}}{\tau^{p^m}}$, $G[i_m(\tau)]=G(u_2)=\gen{\sigma^{p^{n+1}}}{\tau^{p^{m}}}$ and $G[i_{n+1}(\sigma)]=G(u_3)=\gen{\sigma^{p^{n+1}}}{\tau^{p^{m+1}}}$, where $u_1,u_2,u_3$ are three consecutive upper ramification breaks and $G(u)^p=G(u+e)$ for all $u\ge u_1$. Again, let $G'=\langle\sigma,\tau\rangle$ and let $e'$ be the absolute ramification index of $F'$. We have \[\phi_{G'}(i_{n+1}(\sigma))-\phi_{G'}(i_{n}(\sigma))=\frac{i_{m}(\tau)-i_n(\sigma)}{p^{n+m+1}}+\frac{i_{n+1}(\sigma)-i_m(\tau)}{p^{n+m+2}}=e'.\] Similarly, \[\frac{i_{n+1}(\sigma)-i_m(\tau)}{p^{n+m+2}}+\frac{i_{m+1}(\tau)-i_{n+1}(\sigma)}{p^{n+m+3}}=e'.\]
Inductively, we have
\begin{equation}\label{const}
  \frac{i_{m}(\tau)-i_n(\sigma)}{p^{n+m+1}}=\frac{i_{m+j}(\tau)-i_{n+j}(\sigma)}{p^{n+m+1+2j}},\,\forall\,j\in\NN.
\end{equation}
Similarly, we can get \[\frac{i_{n+1}(\sigma)-i_m(\tau)}{p^{n+m+2}}=\frac{i_{n+1+j}(\sigma)-i_{m+j}(\tau)}{p^{n+m+2+2j}},\,\forall\,j\in\NN.\]
This shows \[\frac{i_{n+1}(\sigma)-i_{n}(\sigma)}{i_{n}(\sigma)-i_{n-1}(\sigma)}=p^2\] for all $n$ large enough and we also get the limit $\lim_{n\to\infty}\frac{i_n(\sigma)}{p^{2n}}$ exists.

We summarize this result as the following.

\begin{theorem}\label{height2}
Suppose that $G\subseteq\aut_k(k((x)))$ is a closed subgroup which is isomorphic to $\ZZ_p\times\ZZ_p$. Then the following are equivalent:
\begin{enumerate}
\item For every nonidentity $\sigma\in G$, $\Hgt(\sigma)=2$.
\item For every nonidentity $\sigma\in G$, the sequence $\{\frac{i_n(\sigma)}{p^{2n}}\}_n$ converges.
\item For every nonidentity $\sigma\in G$, $\frac{i_{n+1}(\sigma)-i_n(\sigma)}{i_n(\sigma)-i_{n-1}(\sigma)}=p^2$ for all $n$ sufficiently large.
\item The $\ZZ_p\times\ZZ_p$-extension corresponding to $G$ is of characteristic $0$.
\end{enumerate}
\end{theorem}

\begin{remark}
It is reasonable to extend Theorem \ref{height2} to the case that $G$ is a closed subgroup of $\aut_k(k((x)))$ which is isomorphic to a free $\ZZ_p$-module of rank $n>2$. Our method seems not applicable to show that if every nonidentity element of $G$ has height $n$, then $G$ corresponds to an extension of characteristic $0$. However, in \cite{HsiaLi}, we use different approach to show that the corresponding statements (3) and (4) are equivalent.
\end{remark}

\section{Ramification Index}
Suppose that $G$ is isomorphic to $\ZZ_p$ with generator $\sigma$ and is corresponding to a $\ZZ_p$-extension $E/F$ of characteristic $0$. Then we can use the limit $\lim_{n\to\infty}i_n(\sigma)/p^n$ to determine the absolute ramification index $e$ of $F$. In fact, by (\ref{ram1}) we have $e=\frac{p-1}{p}\lim_{n\to\infty}i_n(\sigma)/p^n$. For the case that $G$ is isomorphic to $\ZZ_p\times \ZZ_p$ with $G=\langle\sigma,\tau\rangle$ and is corresponding to a $\ZZ_p\times\ZZ_p$-extension of characteristic $0$, both limits $\lim_{n\to\infty}i_n(\sigma)/p^{2n}$ and $\lim_{n\to\infty}i_n(\tau)/p^{2n}$ exist, so it is interesting to know whether it is possible to determine the absolute ramification index by purely using the limits $\lim_{n\to\infty}i_n(\sigma)/p^{2n}$ and $\lim_{n\to\infty}i_n(\tau)/p^{2n}$.

For the case of depth $2$, if $G(N)$ is generated by $\sigma^{p^n},\tau^{p^m}$ and $G(u)^p=G(u+e)$ for $u\ge N$, then as indicated above
\[i_{n+j}(\sigma)=i_n(\sigma)+\frac{p^{2+m+n}e}{p^2-1}(p^{2j}-1),\,\forall\,j\in\NN,\]
and hence \[\lim_{j\to\infty}\frac{i_j(\sigma)}{p^{2j}}=\frac{p^2}{p^2-1}p^{m-n}e.\] Similarly, \[\lim_{j\to\infty}\frac{i_j(\tau)}{p^{2j}}=\frac{p^2}{p^2-1}p^{n-m}e.\] Therefore, we have \[e=\frac{p^2-1}{p^2}\sqrt{\lim_{j\to\infty}\frac{i_j(\tau)}{p^{2j}}}\sqrt{\lim_{j\to\infty}\frac{i_j(\sigma)}{p^{2j}}}.\]
Notice that in this case, since $i_{n+j}(\sigma)=i_{m+j}(\tau)$ for all $j\in\NN$, if we set \[\lim_{j\to\infty}\frac{i_j(\sigma)}{p^{2j}}=\gamma_1, \lim_{j\to\infty}\frac{i_j(\sigma)}{p^{2j}}=\gamma_2,\] then $\gamma_1/\gamma_2=p^{2(m-n)}$.
In other words, $\log_p\gamma_1-\log_p\gamma_2$ must be an even number.
\begin{example}
For an odd prime $p$, let $F=\QQ_p(\zeta)$ be the unramified extension of degree $2$ over $\QQ_p$ with $\zeta$ being a unit in $\CO_F$. Consider the Lubin-Tate
formal group over $\CO_F$ constructed by $[p](x)=px+x^{p^2}$. For $\alpha\in\CO_F^*$, let $[\alpha](x)\in\CO_F[[x]]$ be the automorphism of the Lubin-Tate formal group with leading coefficient $\alpha$ and we denote its reduction by $\sigma_\alpha\in\FF_{p^2}[[x]]$. For $\alpha\in \CO_F^*$ with $v_F(\alpha-1)=r$, it is well-known that $i_n(\sigma_\alpha)=p^{2(r+n)}-1$ and hence we have $\lim_{j\to\infty}\frac{i_j(\sigma_\alpha)}{p^{2j}}=p^{2r}$. For the case that $\alpha=1+p$ and $\beta=1+\zeta p$, we know that the closed subgroup generated by $\sigma_\alpha,\sigma_\beta$ corresponds to an extension $N/M$ where $M$ is the extension of $F$ generated by the $p$-torsion elements, i.e elements that satisfy $[p](x)=0$. Therefore, we have the ramification index of $M$ over $\QQ_p$ is $p^2-1$ which is equal to \[\frac{p^2-1}{p^2}\sqrt{\lim_{j\to\infty}\frac{i_j(\sigma_\alpha)}{p^{2j}}}\sqrt{\lim_{j\to\infty}\frac{i_j(\sigma_\beta)}{p^{2j}}}.\] Similarly, for the extension $N/M'$ corresponding the closed subgroup $G'$ generated by $\alpha,\beta^p$, we have the ramification index of $M'$ over $\QQ_p$ is $(p^2-1)p$.
Notice that $G'$ is also generated by $\sigma_{\alpha},\sigma_{\beta'}$ where $\beta'=\alpha\beta^p$, but the ramification index of $M'$ over $\QQ_p$ is not equal to
\[\frac{p^2-1}{p^2}\sqrt{\lim_{j\to\infty}\frac{i_j(\sigma_\alpha)}{p^{2j}}}\sqrt{\lim_{j\to\infty}\frac{i_j(\sigma_{\beta'})}{p^{2j}}}=p^2-1.\] This is because the ramification subgroup of $G(u)$ is not of the form $\gen{\sigma^{p^n}}{\tau^{p^m}}$ when $u$ is large enough.
\end{example}

For the case of depth $1$, if $G(u)$ is generated by $\sigma^{p^n},\tau^{p^m}$ and $G(u')^p=G(u'+e)$ for $u'\ge u$, then as indicated above
\[\frac{i_{m+j}(\tau)-i_{n+j}(\sigma)}{p^{n+m+1+2j}}+\frac{i_{n+1+j}(\sigma)-i_{m+j}(\tau)}{p^{n+m+2+2j}}=e,\,\forall\,j\in\NN.\]
If we set \[\lim_{j\to\infty}\frac{i_j(\sigma)}{p^{2j}}=\gamma_1,\,\lim_{j\to\infty}\frac{i_j(\tau)}{p^{2j}}=\gamma_2,\] then \[e=\frac{p-1}{p^{m-n+1}}\gamma_1+\frac{p-1}{p^{n-m+2}}\gamma_2.\] Moreover, without lose of generality we assume that $i_{n+j}(\sigma)<i_{m+j}(\tau)<i_{n+j+1}(\sigma)$ for all $j\in\NN$. Diving by $p^{2(n+j)}$ and taking limits, we get \[\gamma_1\le p^{2(m-n)}\gamma_2\le p^2\gamma_1.\] However, $(i_{m+j}(\tau)-i_{n+j}(\sigma))/{p^{n+m+1+2j}}$ is a nonzero constant $c$ for all $j\in\NN$ (by (\ref{const})). Taking limits, we get
\[ p^{2(m-n)}\gamma_2-\gamma_1=p^{m-n+1}c\ne 0.\] Similarly, $p^{2(m-n)}\gamma_2\ne p^2\gamma_1$. In other words, $\log_p\gamma_1-\log_p\gamma_2$ cannot be an even number and $m-n$ is the unique integer between $(1/2)(\log_p\gamma_1-\log_p\gamma_2)$ and $1+(1/2)(\log_p\gamma_1-\log_p\gamma_2)$.
\begin{example}
For an odd prime $p$, let $F=\QQ_p(\pi)$ be the totally ramified extension of degree $2$ over $\QQ_p$ with $\pi$ being a prime element in $\CO_F$. Consider the Lubin-Tate
formal group over $\CO_F$ constructed by $[\pi](x)=\pi x+x^{p}$. For $\alpha\in\CO_F^*$, let $[\alpha](x)\in\CO_F[[x]]$ be the automorphism of the Lubin-Tate formal group with leading coefficient $\alpha$ and we denote its reduction by $\sigma_\alpha\in\FF_{p}[[x]]$. For $\alpha\in \CO_F^*$ with $v_F(\alpha-1)=r$, it is well-known that $i_n(\sigma_\alpha)=p^{(r+2n)}-1$ and hence we have $\lim_{j\to\infty}\frac{i_j(\sigma_\alpha)}{p^{2j}}=p^{r}$. For the case that $\alpha=1+\pi$ and $\beta=1+\pi^2$, consider $G$ being the closed subgroup generated by $\sigma_\alpha,\sigma_\beta$. We have \[\gamma_1=\lim_{j\to\infty}\frac{i_j(\sigma_\alpha)}{p^{2j}}=p, \gamma_2=\lim_{j\to\infty}\frac{i_j(\sigma_\beta)}{p^{2j}}=p^2.\] Moreover, $p-1$ is the first upper ramification break of $G$ where $G[p-1]=G(p-1)=G$ is generated by $\sigma_\alpha,\sigma_\beta$ and $G(u)^p=G(u+e)$ for $u\ge p-1$. Notice that $m-n=1-1=0$ is the only integer between $(1/2)(\log_p\gamma_1-\log_p\gamma_2)=-1/2$ and $1+(-1/2)$. On the other hand, $G$
corresponds to an extension $N/M$ where $M$ is the extension of $F$ generated by the $\pi$-torsion elements, i.e elements that satisfy $[\pi](x)=0$. Therefore, we have the ramification index of $M$ over $F$ is $p-1$ and hence the ramification index of $M$ over $\QQ_p$ is $2(p-1)$ which is equal to \[\frac{p-1}{p^{1-1+1}}p+\frac{p-1}{p^{1-1+2}}p^2.\] Similarly, for the extension $N/M'$ corresponding the closed subgroup $G'$ generated by $\alpha,\beta^p$, we have the ramification index of $M'$ over $\QQ_p$ is $2(p-1)p$.
\end{example}

We summarize our result as the following.
\begin{theorem}\label{rami}
Let $G$ be a closed subgroup of $\aut_k(k((x)))$ which is isomorphic to $\ZZ_p\times\ZZ_p$. Suppose $G$ corresponds to an extension of characteristic $0$. Suppose further that $G=\gen{\sigma}{\tau}$, $G(u)=\gen{\sigma^{p^n}}{\tau^{p^m}}$ and $G(u')^p=G(u'+e)$ for all $u'\ge u$.
Let \[\gamma_1=\lim_{j\to\infty}\frac{i_j(\sigma)}{p^{2j}}, \gamma_2=\lim_{j\to\infty}\frac{i_j(\tau)}{p^{2j}}.\]
\begin{enumerate}
\item If $\log_p(\gamma_1/\gamma_2)$ is an even number, then $G$ is of depth $2$ and \[e=\frac{p^2-1}{p^2}\sqrt{\gamma_1\gamma_2}.\]
\item If $\log_p(\gamma_1/\gamma_2)$ is not an even number, then $G$ is of depth $1$. Furthermore, let $a$ be the unique integer between $(1/2)\log_p(\gamma_1/\gamma_2)$ and $1+(1/2)\log_p(\gamma_1/\gamma_2)$. Then \[e=\frac{p-1}{p^{a+1}}\gamma_1+\frac{p-1}{p^{2-a}}\gamma_2.\]
\end{enumerate}
\end{theorem}
\begin{remark}
In the case of depth 2, let $a$ be the integer $(1/2)\log_p(\gamma_1/\gamma_2)$. Then we have
\[\frac{p^2-1}{p^2}\sqrt{\gamma_1\gamma_2}=\frac{p-1}{p^{a+1}}\gamma_1+\frac{p-1}{p^{2-a}}\gamma_2.\]
\end{remark}

\end{document}